\documentclass[12pt,reqno]{amsart}
\usepackage{amssymb,amsmath,amstext,amsthm,amsfonts}
\usepackage[dvips]{graphicx}
\usepackage[ansinew]{inputenc}
\usepackage{a4wide}
\title
{Hyperbolic times: frequency vs. integrability}

\thanks{Work partially supported by ESF through PRODYN
  programme, and FCT through CMUP}



\author{José F. Alves}
\address{Centro de Matem\'atica da Universidade do
Porto\newline \indent Rua do Campo Alegre 687, 4169-007 Porto, Portugal} \email{jfalves@fc.up.pt}
\urladdr{www.fc.up.pt/cmup/jfalves}

\author{V{\'\i}tor Ara\'ujo}
\address{Centro de Matem\'atica da Universidade do Porto
\newline \indent Rua do Campo Alegre 687, 4169-007 Porto, Portugal}
\email{vdaraujo@fc.up.pt} \urladdr{www.fc.up.pt/cmup/vdaraujo}

\date{\today}

\begin{document}

\newtheorem{theorem}{Theorem}
\newtheorem{corollary}{Corollary}

\newcommand{\mcup}{\mbox{$\bigcup$}}
\newcommand{\mcap}{\mbox{$\bigcap$}}

\newcommand{\leb}{\operatorname{Leb}}
\newcommand{\dist}{\operatorname{dist}}

\def \RR {{\mathbb R}}
\def \ZZ {{\mathbb Z}}
\def \NN {{\mathbb N}}
\def \PP {{\mathbb P}}
\def \TT {{\mathbb T}}

\def \ra {\rightarrow }
 \def \wh {\widehat }
 \def \un{\underline }
 \def \ov {\overline}
 \def \supp {\mbox{supp}\, }
 \def \wlim {\mbox{$w^*$-}\lim_{n\ra\infty}\, }
 \def \distp{\mbox{d}_\PP }

\def \al {\alpha } \def \be {\beta } \def \de {\delta }
\def \ga {\gamma } \def \vare {\varepsilon } \def \vfi {\varphi
} \def \th {\theta } \def \si {\sigma } \def \ep {\varepsilon}

 \def \cf {\mathcal{F}}
 \def \cm {\mathcal{M}}
 \def \cn {\mathcal{N}}
 \def \cq {\mathcal{Q}}
 \def \cp {\mathcal{P}}
 \def \cc {\mathcal{S}}
 \def \ch {\mathcal{H}}
 \def \cU {\mathcal{U}}
  \def \cs {\mathcal{S}}
  \def \cO {\mathcal{O}}

\newcommand{\dem}{\begin{proof}}
\newcommand{\cqd}{\end{proof}}

\newcommand{\qand}{\quad\text{and}\quad}

\newtheorem{maintheorem}{Theorem}
\renewcommand{\themaintheorem}{\Alph{maintheorem}}
\newcommand{\cmt}{\begin{maintheorem}}
\newcommand{\fmt}{\end{maintheorem}}

\newtheorem{maincorollary}[maintheorem]{Corollary}
\newcommand{\cmc}{\begin{maincorollary}}
\newcommand{\fmc}{\end{maincorollary}}

\newtheorem{T}{Theorem}[section]
\newcommand{\ct}{\begin{T}}
\newcommand{\ft}{\end{T}}

\newtheorem{Corollary}[T]{Corollary}
\newcommand{\cco}{\begin{Corollary}}
\newcommand{\fco}{\end{Corollary}}

\newtheorem{Proposition}[T]{Proposition}
\newcommand{\cpr}{\begin{Proposition}}
\newcommand{\fpr}{\end{Proposition}}

\newtheorem{Lemma}[T]{Lemma}
\newcommand{\cle}{\begin{Lemma}}
\newcommand{\fle}{\end{Lemma}}

\theoremstyle{remark}

\newtheorem{Remark}[T]{Remark}
\newcommand{\cre}{\begin{Remark}}
\newcommand{\fre}{\end{Remark}}

\newtheorem{Definition}{Definition}
\newcommand{\cde}{\begin{Definition}}
\newcommand{\fde}{\end{Definition}}

\maketitle

\begin{abstract}

We consider  dynamical systems on compact manifolds, which are
local diffeomorphisms outside an exceptional set (a compact
submanifold). We are interested in analyzing the relation between
the integrability (with respect to Lebesgue measure) of the first
hyperbolic time map  and the existence of  positive frequency of
hyperbolic times. We show that some (strong) integrability of the
first hyperbolic time map implies   positive frequency of
hyperbolic times. We also present an example of a map with
positive frequency of hyperbolic times at Lebesgue almost every
point but whose first hyperbolic time map is not integrable with
respect to the Lebesgue measure.
\end{abstract}

\tableofcontents

\section{Introduction}
\label{sec.intro}

In the last decades many dynamicists have dedicated their
attention  to the understanding of the dynamical features of
systems exhibiting some non-uniformly hyperbolic behavior. 
In this direction we mention \cite{BC1, BY1,J} for quadratic maps, \cite{BC2,BY2,MV, V1,WY} for Hénon-like
diffeomorphisms, and \cite{Al,AV,V} for a generalized higher dimensional version of the quadratic and Hénon-like
maps. The dynamics of all these systems is characterized by the existence of regions of the phase space where
the system displays some hyperbolicity, together with {\em critical regions} where some strong non-hyperbolic
behavior appears. The strategy for dealing with the loss of hyperbolicity on the quadratic and Hénon-like maps
is based on  the existence of well-defined {\em recovering periods} during which the non-hyperbolic effect of
the critical region is compensated for.

The recovering period argument is no longer possible to be
used in the class of endomorphisms introduced in \cite{V}
due to the fact that positive iterates of the critical
region have unavoidable intersections with this region.
Instead, the mechanism that enables one to obtain the
non-uniformly expanding behavior in that case is of a
statistical type and comes from a delicate analysis on the
derivative of the map along full orbits.

A new and powerful tool in this setting has been introduced in \cite{Al} through the notion of {\em hyperbolic
times}. These have become a very useful ingredient in the study of non-hyperbolic dynamical systems, playing an
important role  in the proof of several results about the existence of absolutely continuous invariant measures
and their statistical properties; see \cite{Al,AA,ABV,ALP,AV,BBM}. Ideas of hyperbolic times were implicitly
contained in Pesin's theory and in the work of Pliss and Mañé. Also, recently hyperbolic times have been used to
study stochastic flows in \cite{DKK}.

The applications of hyperbolic times have been twofold: on the one hand, some of these works deal with the
integrability with respect to the Lebesgue measure of the first hyperbolic time map, while on the other hand,
the usefulness of hyperbolic times appears through their positive frequency along typical orbits.

Our aim in this work is to clarify the relations among the
integrability of the first hyperbolic time map, the frequency of
hyperbolic times and the existence of absolutely continuous
invariant measures.

\subsection*{Statement of results}
We start by introducing the most relevant concepts and
definitions.  Let $f:M\ra M$ be a continuous map defined on
a compact Riemannian manifold $M$ with the induced distance
that we denote by $\dist$, and fix a normalized Riemannian
volume form $m$ on $M$ that we call {\em Lebesgue measure}.
%
%

 Throughout this work we will assume that \( f \) is a local
 diffeomorphism in all of $M$ but an exceptional set
 \( \cs \subset M\), where $\cs$ is a compact submanifold of
 $M$ with $\dim(\cs)<\dim(M)$ (thus $m(\cs)=0$) satisfying
some non-degeneracy conditions. Examples of systems satisfying Definition~\ref{d.nd} include one-dimensional
quadratic maps and
 Viana maps \cite{V}.

\cde\label{d.nd} We say that \( \cs \subset M\) is a {\em
non-degenerate exceptional set} for $f$ if  the following
conditions hold. The first one essentially says that
 $f$  {\em behaves like a power of the distance}
 to \( \cs \):
 there are constants $B>1$ and $\be>0$ such that for every $x\in
 M\setminus\cs$
\begin{enumerate}
 \item[(s$_1$)]
\quad $\displaystyle{\frac{1}{B}\dist(x,\cs)^{\be}\leq \frac
{\|Df(x)v\|}{\|v\|}\leq B\dist(x,\cs)^{-\be}}$ for all $v\in T_x
M$.
\end{enumerate}
Moreover, we assume that the functions \(  \log|\det Df(x)| \) and
\( \log \|Df(x)^{-1}\| \) are \emph{locally Lipschitz} at points
\( x\in M \setminus \mathcal S \) with Lipschitz constant
depending on $ \dist (x, \mathcal S)\): for every $x,y\in
M\setminus \cs$ with $\dist(x,y)<\dist(x,\cs)/2$ we have
\begin{enumerate}
\item[(s$_2$)] \quad $\displaystyle{\left|\log\|Df(x)^{-1}\|-
\log\|Df(y)^{-1}\|\:\right|\leq
\frac{B}{\dist(x,\cs)^{\be}}\dist(x,y)}$;
 \item[(s$_3$)]
\quad $\displaystyle{\left|\log|\det Df(x)|- \log|\det
Df(y)|\:\right|\leq \frac{B}{\dist(x,\cs)^{\be}}\dist(x,y)}$.
 \end{enumerate}
 \fde
The set $\cs$ may be taken as some set of critical points of $f$
or a set where $f$
 fails to be differentiable. The case where
 $\cs$ is equal to the empty set may also be considered.
For the next definition it will be useful to introduce \(
\dist_{\delta}(x,\cs) \), the \( \delta \)-\emph{truncated}
distance from \( x \) to~$\cs$, defined as \(
\dist_{\delta}(x,\cs) =  \dist(x,\cs) \) if \( \dist(x,\cs) \leq
\delta\), and \( \dist_{\delta}(x,\cs) =1 \) otherwise.


\cde Let  $\beta>0$ be as in Definition~\ref{d.nd}, and fix  $b>0$   such that $b < \min\{1/2,1/(4\beta)\}$.
Given $0<\sigma<1$ and $\delta>0$, we say that $n$ is a {\em $(\sigma,\delta)$-hyperbolic time\footnote{In the
case $\cs=\emptyset$ the definition of $(\sigma,\delta)$-hyperbolic time reduces to the first condition in
\eqref{d.ht}, and we simply call it $\sigma$-hyperbolic time. }} for a point $x\in M$ if for all $1\le k \le n$,
 \begin{equation}\label{d.ht}
\prod_{j=n-k}^{n-1}\|Df(f^j(x))^{-1}\| \le \sigma^k \qand
\dist_\delta(f^{n-k}(x), \cs)\ge \sigma^{b k}.
 \end{equation}
 We say that the {\em frequency of
$(\sigma,\delta)$-hyperbolic times} for $x\in M$ is greater than
$\theta>0$ if, for large $n$, there are $n_1<n_2\dots <n_\ell\le
n$ which are $(\sigma,\delta)$--hyperbolic times
 for $x$ and $\ell \ge\th n$. 
\fde

We point out that  condition \eqref{d.ht} implies the dynamically meaningful property
$$\|(Df^{k}(f^{n-k}(x)))^{-1}\| \leq \sigma^k,\quad\text{for $1\le k\le n$},$$ which says that at any intermediate moment of time between 0
and $n$ we get exponential stretching by iteration under $f$. The work of Viana \cite{V} provides interesting
higher dimensional examples of maps with many hyperbolic times for most points. The significance of hyperbolic
times may be attested by the following result whose proof is essentially contained in~\cite{ABV}.

\cmt\label{t.abv} Let $f\colon M\to M$ be a $C^2$ local
diffeomorphism outside a non-degenerate exceptional set
$\cs\subset M$. If there are $0<\sigma<1$,
   $\delta>0$,  and $H\subset M$ with $m(H)>0$ such that
the frequency of $(\sigma,\delta)$-hyperbolic  times is bigger
than $\theta>0$ for every $x\in H$, then $f$ has some absolutely
continuous invariant probability measure.
 \fmt

The existence of $(\sigma,\delta)$-hyperbolic times for  Lebesgue almost all points in $M$ allows us to
introduce a map $h\colon M\to \ZZ^+$ defined Lebesgue almost everywhere and assigning to each $x\in M$ its first
$(\sigma,\delta)$-hyperbolic time. The integrability properties of this first hyperbolic time map play an
important role in the study of some statistical properties of several classes of dynamical systems, such as
stochastic stability and decay of correlations; see \cite{AA,ALP,AV,BY1}. The same conclusion of
Theorem~\ref{t.abv} can be obtained under the assumption of integrability with respect to Lebesgue measure of
the first hyperbolic time map.

 \cmt\label{t.integrable-h-1} Let $f\colon M\to M$ be a $C^2$ local
diffeomorphism outside a non-degenerate exceptional set $\cs\subset M$.  If for some $0<\sigma<1$ and
   $\delta>0$, the first $(\sigma,\delta)$-hyperbolic time map $h:M\to\ZZ^+$ is
Lebesgue integrable, then $f$ has
  an  absolutely continuous invariant
probability measure $\mu$.
 \fmt

 Having in mind Theorem~\ref{t.abv} and Theorem~\ref{t.integrable-h-1},  one
is naturally  interested in understanding
 the relation between the
existence of a positive Lebesgue measure subset of points in
 $M$
with positive frequency of $(\sigma,\delta)$-hyperbolic times and
the integrability with respect
 to the Lebesgue measure of the first $(\sigma,\delta)$-hyperbolic time map.

 \cmt\label{t.integrable-h-2}  Let $f\colon M\to M$ be a $C^2$ local
diffeomorphism outside a non-degenerate critical set $\cs\subset
M$. If for $0<\sigma<1$ and
   $\delta>0$ the first $(\sigma,\delta)$-hyperbolic time map $h:M\to\ZZ^+$  belongs to
$L^p(m)$ for some  $p>4$, then there are $\hat\sigma>0$ and
$\theta>0$ such that Lebesgue almost every $x\in M$ has frequency
of $(\hat\sigma,\delta)$-hyperbolic times bigger than $\theta$.
 \fmt

 We do not know if the need for stronger integrability in
 this last theorem is due to the methods we have used to
 prove it or some kind of stronger integrability is really
 necessary. It remains an interesting open question to know
 the smallest value of $p\ge 1$ for which the first
 condition in Theorem~\ref{t.integrable-h-2} still implies
 the desired conclusion. As a by-product of the proof of
 Theorem~\ref{t.integrable-h-2} we will see that the answer
 to this question is optimal when $\cs=\emptyset$.

\cmt\label{co:t-integrable-h} Let $f\colon M\to M$ be a $C^2$ local diffeomorphism. If for some $\sigma\in(0,1)$
the first $\sigma$-hyperbolic time map is Lebesgue integrable, then there are $\hat\sigma>0$ and $\th>0$   such
that Lebesgue almost every $x\in M$ has frequency of $\hat\sigma$-hyperbolic times bigger than $\theta$.  \fmt

In the opposite direction, one could ask whether the positive
frequency of hyperbolic times is enough for assuring the
integrability of the first hyperbolic time map. There is no hope
of such a result. In Section~\ref{sec.example} we present an
example of a map of the circle (with nonempty exceptional set)
having positive frequency of hyperbolic times for Lebesgue almost
every point, but whose first hyperbolic time map is not integrable
with respect to the Lebesgue measure.

Another interesting open question is whether a $C^2$ local
diffeomorphism (with no exceptional set) on a compact manifold,
admitting positive frequency of hyperbolic times for Lebesgue
almost every point, necessarily has a first hyperbolic time map
that is Lebesgue integrable.

Hyperbolic times appear naturally when $f$ is assumed to be {\em
non-uniformly expanding} in some set $H\subset M$: there is some
$c>0$ such that for every $x\in H$ one has
 \begin{equation} \label{liminf1}
\limsup_{n\ra +\infty}\frac{1}{n}
\sum_{j=0}^{n-1}\log\|Df(f^j(x))^{-1}\|<-c,
 \end{equation}
and points in $H$ satisfy some {\em slow recurrence to the
exceptional set}: given any $\vare>0$ there exists $\delta>0$ such
that for every $x\in H$
\begin{equation} \label{e.faraway1}
    \limsup_{n\to+\infty}
\frac{1}{n} \sum_{j=0}^{n-1}-\log \dist_\delta(f^j(x),\cs)
\le\vare.
\end{equation}
The next result has been  proved in \cite{ABV} (see Theorem~C and
Lemma 5.4 therein) and will be used
 in the proof  of our results.

\cmt\label{t.abv2} Let $f\colon M\to M$ be a $C^2$ local diffeomorphism outside a non-degenerate exceptional set
$\cs\subset M$. If  there is some set $H\subset M$ with $m(H)>0$ such that \eqref{liminf1} and
\eqref{e.faraway1} hold for all $x\in H$, then there are $0<\sigma<1$, $\delta>0$ and $\theta>0$ such that the
frequency of $(\sigma,\delta)$-hyperbolic times for the points in $H$ is bigger than $\theta$.
 \fmt

This work is organized as follows. In Section~\ref{sec.positive} we establish the basic properties of hyperbolic
times and sketch the proof of Theorem~\ref{t.abv}, using some of the results in \cite{ABV}. In
Section~\ref{sec:integrability} we prove Theorem~\ref{t.integrable-h-1}, and in Section~\ref{sec.intfreq} we
prove Theorem~\ref{t.integrable-h-2} and Theorem~\ref{co:t-integrable-h}. In Section~\ref{sec.example} we
present an example of a map of the circle which is a $C^2$ local diffeomorphism everywhere  but in an
exceptional set with two points, having positive frequency of hyperbolic times, and whose first hyperbolic time
map is non-integrable with respect to the Lebesgue measure.

\medskip
\noindent\emph{Acknowledgments.}  We are indebted to Henk Bruin,
who suggested us the example of Section~\ref{sec.example}  in
conversations during the workshop \emph{Concepts and Techniques in
Smooth Ergodic
 Theory} held at the Imperial College, London,
 July 2001. We also thank Xavier Bressaud and Sandro Vaienti for
 helpful conversations at the Institut de Mathématiques de Luminy,
 Marseille, June 2002.


\section{Positive frequency of hyperbolic times}
\label{sec.positive}

One of the main features of hyperbolic times is that the
corresponding iterates locally behave as those of an expanding
map, namely with uniform expansion and uniformly bounded
distortion. This is precisely stated in the next result.

 \cle\label{p.contr}
 There are $\delta_1>0$ and $C_1>0$
 such that if $n$ is a
$(\sigma,\delta)$-hyperbolic time for $x$, then there is a
neighborhood $V_x$ of $x$ in $M$ for which
 \begin{enumerate}
 \item $f^n$ maps $V_x$ diffeomorphically onto
 the ball of radius $\delta_1$ around $f^n(x)$;
 \item for $1\le k <n$ and $y,
z\in V_x$, $ \dist(f^{n-k}(y),f^{n-k}(z)) \le
\sigma^{k/2}\dist(f^{n}(y),f^{n}(z))$;
 \item $f^n\vert V_x$ has   distortion bounded by $C_1$: if  $y, z\in
V_x$, then
 $$ \frac{1}{C_1} \le \frac{|\det Df^n (y)|}{|\det
Df^n (z)|}\le C_1 \,. $$
  \end{enumerate}
 \fle
 \dem See \cite[Lemma 5.2]{ABV} and \cite[Corollary
5.3]{ABV}. \cqd

Let us now give a brief idea on the way we  obtain Theorem~\ref{t.abv} from the results  in \cite{ABV}.

Assume that there are $0<\sigma<1$,
   $0<\delta$,  and $H\subset M$ with $m(H)>0$ such that
for every $x\in H$ the frequency of $(\sigma,\delta)$-hyperbolic
times is bigger than $\theta>0$. Given an integer $n \ge 1$ we
define
 $$
 H_n=\{ x\in H\colon \mbox{ $n$ is a
 $(\sigma,\delta)$-hyperbolic time for $x$}\}.
 $$

\cpr \label{p.maximal} Take $\delta_1$  as in Lemma~\ref{p.contr}. There exists a constant $\tau>0$ such that
for any $n$ there exists a finite subset $\widehat{H}_n$ of $H_n$ for which the balls of radius $\delta_1/4$
around the points $x\in f^n(\widehat{H}_n)$ are pairwise disjoint, and their union $\Delta_n$ satisfies
$$
f_*^n m(\Delta_n \cap H) \ge f_*^n m(\Delta_n \cap f^n(H_n)) \ge
\tau m(H_n)
$$
\fpr \dem See \cite[Proposition 3.3]{ABV}. \cqd

 From
Proposition~\ref{p.maximal}, we may find for each $j\ge 1$ a
finite set of points $x_1^j, \ldots, x_N^j$ (in principle with $N$
depending on $j$) admitting $j$ as a $(\sigma,\delta)$-hyperbolic
time, such that:
\begin{enumerate}
\item $V_{x_1^j}, \ldots, V_{x_N^j}$ are pairwise disjoint; \item the Lebesgue measure of $W_j=V_{x_1^j} \cup
\ldots \cup V_{x_N^j}$ is larger than the Lebesgue measure of $H_j$, up to a uniform multiplicative constant
$\tau>0$.
\end{enumerate}
We let $(\mu_n)_n$ be the sequence of the averages of the positive
iterates of Lebesgue measure on $M$,
\begin{equation*}
 \mu_n=\frac{1}{n} \sum_{j=0}^{n-1} f_*^j m,
 \end{equation*}
and $\nu_n$ be the part of $\mu_n$ carried on disks of radius
$\delta_1$ around points $f^j(x_k^j)$ such that $1\le j\le n$ is a
$(\sigma,\delta)$-hyperbolic time for $x$,
$$
\nu_n=\frac{1}{n}\sum_{j=0}^{n-1}f_*^j(m \mid W_j).
$$
By Proposition~\ref{p.maximal} we have
$$\nu_n(H)\ge\frac\tau n\sum_{i=0}^{n-1}m(H_i).$$ So, it suffices to prove
that this last expression is larger than some positive constant,
for $n$ large.
 Let $\xi_n$ be the measure in
$\{1,\dots,n\}$ defined by $\xi_n(B)=\# B/n$, for each subset $B$.
Then, using Fubini's theorem
\begin{eqnarray*}
\frac{1}{n} \sum_{i=0}^{n-1}m(H_n)
& =& \int \left(\int \chi(x,i)\,dm(x)\right)d\xi_n(i) \\
& = & \int \left(\int \chi(x,i)\,d\xi_n(i)\right)dm(x),
\end{eqnarray*}
where $\chi(x,i)=1$ if $x\in H_i$ and $\chi(x,i)=0$ otherwise.
Now, since we are assuming positive frequency of hyperbolic times
for points in $H$, this means that the integral with respect to
$d\xi_n$ is larger than $\theta>0$ for large $n$. So, the
expression on the right hand side is bounded from below by $\theta
m(M)$. This implies that each $\nu_n$ has total mass uniformly
bounded away from zero. Moreover, as a consequence of the bounded
distortion given by Lemma~\ref{p.contr}, every $f_*^j(m \mid W_j)$
is absolutely continuous with respect to Lebesgue measure, with
density uniformly bounded from above, and so the same is true for
every $\nu_n$.

Since we are working with a continuous map in the compact space
$M$, we know that sequences of probability measures in $M$ have
weak* accumulation points.  Take $n_k\to\infty$ such that both
$\mu_{n_k}$ and $\nu_{n_k}$ converge in the weak$^*$ sense to
measures $\mu$ and $\nu$, respectively. Then $\mu$ is an invariant
probability measure, $\mu=\nu+\eta$ for some measure $\eta$, $\nu$
is absolutely continuous with respect to Lebesgue measure, and
$\nu(H)>0$. Now, if $\eta=\eta_{ac}+\eta_{s}$ denotes the Lebesgue
decomposition of $\eta$ (as the sum of an absolutely continuous
and a completely singular measure, with respect to Lebesgue
measure), then $\mu_{ac}=\nu+\eta_{ac}$ gives the absolutely
continuous component in the corresponding decomposition of $\mu$.
By uniqueness of the Lebesgue decomposition, and the fact that the
push-forward under $f$ preserves the class of absolutely
continuous measures, we may conclude that $\mu_{ac}$ is an
invariant measure. Clearly, $\mu_{ac}(H)\ge \nu(H)>0$. Normalizing
$\mu_{ac}$ we obtain an absolutely continuous $f$-invariant
probability measure.

The next lemma will be useful in Section~\ref{sec.example}.

\cle \label{l.disco} Take $\delta_1$  as in Lemma~\ref{p.contr}. If $H\subset M$ is a positively invariant set
with $m(H)>0$ for which \eqref{liminf1} and \eqref{e.faraway1} hold, then there exists some disk $\Delta$ with
radius $\delta_1/4$ such that $m(\Delta \setminus H)=0$. \fle \dem See \cite[Lemma 5.6]{ABV}.\cqd

Using this lemma it is shown in \cite[Section 5]{ABV}  that $f$ has a finite number of absolutely continuous
invariant probability
measures  which are ergodic. 

\section{Integrability of first hyperbolic time map}
\label{sec:integrability}

Here we prove Theorem~\ref{t.integrable-h-1}. As in
Section~\ref{sec.positive} the strategy is to consider $(\mu_n)_n$
the sequence of averages of forward iterates of Lebesgue measure
on $M$
\begin{equation*}
 \mu_n=\frac{1}{n} \sum_{j=0}^{n-1} f_*^j m.
 \end{equation*}
 Since we are dealing with a continuous map of a compact manifold,
 we know that the sequence $(\mu_n)_n$ has accumulation points -- which belong to the space of
 $f$-invariant probability measures.
 Now the idea is to show that  such accumulation points are absolutely continuous
 with respect to the Lebesgue measure. 

\cpr\label{p.density}
 There is a constant $C_2>0$ (depending on $\delta_1$ and $C_1$ from Lemma~\ref{p.contr})
  such that for
every $n\geq 0$
 $$
 \frac{d}{dm}f^n_*\big(m\mid H_n\big)\leq C_2.
 $$
 \fpr
 \dem Take $\delta_1>0$ given by Lemma~\ref{p.contr}.
It suffices to show that there is some uniform constant $C>0$ such
that if  $A\subset M$ is a Borel set with diameter smaller than
$\delta_1/2$, then
 $$
 m\big(f^{-n}(A)\cap  H_n\big)\leq C m(A).
 $$
 Let $A$ be a Borel set in $M$ with diameter smaller than
$\delta_1/2$ and $B$ an open  ball of radius $\delta_1/2$
containing $A$. Taking the connected components of $f^{-n}(B)$ we
may write
 $$
 f^{-n}(B)=\bigcup_{k\geq 1}B_k,
 $$
where $(B_k)_{k\geq 1}$ is a (possibly finite) family of pairwise disjoint open sets in $M$. Taking into account
only those $B_k$ that  intersect $H_n$, we choose, for each $k\geq 1$, a point $x_k\in H_n\cap B_k$. For each
$k\geq 1$ let $V_{x_k}$ be the neighborhood of $x_k$  given by Lemma~\ref{p.contr}. Since $B$ is contained in
$B\big(f^n(x_k), {\delta_1}\big)$, the ball of radius $\delta_1$ around $f^n(x_k)$, and $f^n$ is a
diffeomorphism from $V_{x_k}$ onto $B\big(f^n(x_k), {\delta_1}\big)$, we must have $B_k\subset V_{n_k}$ (recall
that by the choice of $B_k$ we have $f^n(B_k)\subset B$). As a consequence of this and Lemma~\ref{p.contr} we
have that $f^n\mid B_k\colon B_k\rightarrow B$ is a diffeomorphism with uniform bounded distortion for all
$n\ge1$ and $k\ge1$:
 $$ \frac{1}{C_1}\leq \frac{|\det Df^n(y)|}{|\det Df^n(z)|}
  \leq C_1\quad\text{for all $y,z\in B_k$.} $$
 This finally
gives
 \begin{eqnarray*}
  m\big(f^{-n}(A)\cap H_n\big) &\leq &
  \sum_{k\ge1}m\big(f^{-n}(A\cap B)\cap B_k\big)\\
   &\leq & \sum_{k\ge1}C_1\frac{m(A\cap B)}{m(B)}m(B_k)\\
    &\leq & C_2 m(A),
 \end{eqnarray*}
for some constant $C_2>0$ only depending on $C_1>0$ and on the
volume of the ball $B$ of radius $\delta_1/2$. \cqd

Defining, for each $n\ge 1$,
$$
 H_n^*=\{ x\in M\colon \mbox{ $n$ is the {\em first}
 $(\sigma,\delta)$-hyperbolic time for $x$}\},
 $$
we immediately have
 \begin{equation}\label{eq.unif}
 \int_Mhdm=\sum_{k=1}^\infty km(H^*_k).
 \end{equation}
It will be useful to  define, for each $n,k\ge1$,
 $$
 R_{n,k}= \big\{\,x\in H_n\:\colon \:f^n(x)\in H^*_k\: \big\}.
 $$
Observe that $R_{n,k}$ is precisely the set of  points $x\in M$
for which $n$ is a $(\sigma,\delta)$-hyperbolic time and $n+k$ is
the next $(\sigma,\delta)$-hyperbolic time for $x$ after $n$.
Defining the measures
 \begin{equation}\label{e.nu}
 \nu_n= f^n_*(m\mid H_n)
 \end{equation}
 and
 \begin{equation}\label{e.eta}
\eta_n=\sum_{k=2}^\infty\sum_{j=1}^{k-1} f^{n+j}_*(m\mid R_{n,k}),
 \end{equation}
we may write
 $$
 \mu_n\leq
\frac{1}{n}\sum_{j=0}^{n-1}(\nu_j+\eta_j).
 $$
It follows from  Propositions~\ref{p.density} that
 \begin{equation}\label{dens}
 \frac{d\nu_n}{dm}\leq C_2
 \end{equation}
 for every $n\geq 0$, with $C_2$  not depending on  $n$.
  Our goal now is to
 control the densities of the measures $\eta_n$.

 \cpr\label{p.dens}
 Given $\vare>0$, there is
$C_3(\vare)>0$ such that for  every $n\geq 1$ we may bound
$\eta_n$ by the sum of two non-negative measures, $\eta_n \leq
\omega+\rho$, with
 $$
 \frac{d\omega}{dm}\leq C_3(\vare)\quad\mbox{and}\quad
 \rho(M)<\vare.
 $$
 \fpr
\dem  Let $A$ be some Borel set in $M$.
 For each $n\geq 0$ we have
 \begin{eqnarray*}
 \eta_n(A)
 &=&
 \sum_{k=2}^\infty\sum_{j=1}^{k-1}
 m\big( f^{-n-j}(A)\cap R_{n,k}\big)\\
 &\leq &
 \sum_{k=2}^\infty\sum_{j=1}^{k-1}m\big(f^{-n}\big(f^{-j}(A)
 \cap H^*_k\big)\cap
 H_n\big)\\
 &\leq &
 \sum_{k=2}^\infty\sum_{j=1}^{k-1} C_2
 m\big( f^{-j}(A)\cap
 H^*_k\big).
 \end{eqnarray*}
 (in this last inequality we have used the bound (\ref{dens}) above). Let
now $\vare>0$ be some fixed small number. By the integrability of
$h$ and  since (\ref{eq.unif}) holds, we may choose some  integer
$\ell=\ell(\vare)$ for which
 $$
 \sum_{j=\ell}^{\infty}k\,
 m\big(H^*_k\big)<\frac{\vare}{C_2}.
 $$
We take
 $$
 \omega=C_2\sum_{k=2}^{\ell-1}\sum_{j=1}^{k-1}
 f^j_*(m\mid H^*_k)
 $$
  and
 \begin{equation}\label{e.ro}
 \rho=C_2\sum_{k=\ell}^\infty\sum_{j=1}^{k-1}
 f^j_*(m\mid H^*_k).
 \end{equation}
This last measure satisfies
 $$
 \rho(M)=C_2\sum_{k=\ell}^\infty\sum_{j=1}^{k-1}
 m( H^*_k)\leq C_2
 \sum_{k=\ell}^\infty k \, m( H^*_k)
 <\vare.
 $$
On the other hand,  we have
 $$
 \omega\leq C_2\sum_{k=2}^{\ell-1}\sum_{j=1}^{k-1}
 f^j_*m,
 $$
 and this last measure has density bounded by some
 constant by the non-degeneracy conditions of $f$, since we
 are taking a finite number of push-forwards of Lebesgue
 measure.  \cqd

 It follows from this last proposition and (\ref{dens})
 that weak$^*$
 accumulation points of $(\mu_n)_n$ cannot
 have singular part, thus being absolutely continuous with
 respect to the Lebesgue measure. Since such weak$^*$
 accumulation points are invariant with respect to $f$, we have
 proved Theorem~\ref{t.integrable-h-1}.

 \cre
This argument proves that \emph{every} weak* accumulation point of $(\mu_n)_n$ is absolutely continuous with
respect to Lebesgue measure, whenever the first hyperbolic time function is integrable. This is not known if we
only assume positive frequency of hyperbolic times.
 \fre


\section{Strong integrability implies positive frequency}
\label{sec.intfreq}

Assume that $f:M\to M$ is a $C^2$ local diffeomorphism outside a
non-degenerate exceptional set $\cs\subset M$. We start the proof
of Theorem~\ref{t.integrable-h-2} by obtaining a simple useful
result.

\cpr\label{le:logdist-int} If $\cs$ is a compact submanifold of
$M$ with $\dim(\cs)<\dim(M)$, then the function $\log\dist(x,\cs)$
belongs to $L^p(m)$ for every $1\le p<\infty$. \fpr

\dem We may assume without loss of generality that $\cs$ is
connected. Let $\dim(\cs)=k< n=\dim(M)$. We may cover $\cs$ with
finitely many images of charts $\psi_i(U_i)$ ($i=1,\dots,p$) such
that $U_i\subset\RR^n$ is a bounded open set and
$\psi_i^{-1}(\cs)\subset U_i\cap(\RR^k\times 0^{n-k})$. Denoting
by $\lambda$ the usual $n$-dimensional volume on $\RR^n$ and by
$d$ the standard Euclidean distance on $\RR^n$, then there are
constants $C,K>0$ such that for all $i=1,\dots,p$
\[
\frac1C\le\frac{d(\psi_i^{-1})_* m}{d\lambda}\le C,
\]
and  for all  $w,z\in U_i$
 $$
 \frac1K d(w,z)\le
{\dist(\psi_i(w),\psi_i(z))}\le K d(w,z).
 $$
Hence, for showing that $\log\dist(x,\cs)$ is integrable with
respect to $m$, it is enough to show that $\log
d(x,U\cap(\RR^k\times 0^{n-k}))$ is integrable with respect to
$\lambda$ for any open and bounded neighborhood $U$ of the origin
in $\RR^n$. We may assume without loss of generality that $U$ is
sufficiently small in order to $U\subset B_k\times B_{n-k}$, where
$B_k$ and $B_{n-k}$ are the unit balls around the origin in
$\RR^k$ and $\RR^{n-k}$ respectively. For
$z=(z_1,\dots,z_n)\in\RR^n$ we have
$$d(z, \RR^k\times 0^{n-k})=(z_{k+1}^2+\dots+z_n^2)^{1/2}.$$
Hence, we have for $1\le p<\infty$
\[
\int_{U} |\log d(z,\RR^k\times 0^{n-k})|^p d\lambda \le
\frac1{2^p} \int_{B_{k}}\left(
\int_{B_{n-k}}\!\!\!|\log(z_{k+1}^2+\dots+z_n^2)|^p dz_{k+1}\cdots
dz_n\right) dz_1\cdots dz_k.
\]
 Now it is enough to show that the
inner integral  in the last expression is finite. Actually,
denoting by $S^{n-k-1}_\rho$ the $(n-k-1)$-sphere with radius
$\rho$ around the origin in $\RR^{n-k}$, $dA$ its area element and
$a$ the total area of $S^{n-k-1}_1$, we have
\begin{eqnarray*}
\int_{B_{n-k}}\!\!\!|\log(z_{k+1}^2+\dots+z_n^2)|^p dz_{k+1}\cdots
dz_n & = & \int_0^1 \left( \int_{S^{n-k-1}_\rho} |2 \log \rho|^p
\, dA \right) \, d\rho
\\
& = & a\int_0^1 \rho^{n-k-1}|\log\rho|^p \, d\rho .
\end{eqnarray*}
Since this last integral is finite, we have completed the proof of
the result. \cqd

Assume now that $h$ is integrable with respect to the Lebesgue
measure. By Theorem~\ref{t.integrable-h-1} we know that there
exists an absolutely continuous invariant probability measure
$\mu$ for $f$.

\cco\label{c.integra} If the density $d\mu/dm$ belongs to $L^q(m)$
for some $q>1$, then $\log\dist(x,\cs)$ is $\mu$-integrable. \fco

\dem This is an immediate application of H\"older inequality.
Actually, since
 $$
\int \log\dist(x,\cs)d\mu = \int
\log\dist(x,\cs)\frac{d\mu}{dm}\,dm,
 $$
and we have $d\mu/dm$ in $L^q(m)$ for some $q>1$ and
$\log\dist(x,\cs)$ in $L^p(m)$ for every $p$, then taking $p$
equal to the conjugate of $q$, that is $p^{-1}+q^{-1}=1$, then
H\"older inequality gives that the integral above is finite.
 \cqd

We are also interested in obtaining the same conclusion of the
previous corollary under the hypothesis that $h\in L^p(m)$ for
some $p>4$. Observe that the absolutely continuous $f$-invariant
measure $\mu$ may be obtained as a weak* accumulation point of the
sequence $(\mu_n)_n$ of averages of push-forwards of Lebesgue
measure. As shown in Section~\ref{sec:integrability}, we may write
$$
 \mu_n\leq
\frac{1}{n}\sum_{j=0}^{n-1}(\nu_j+\eta_j),
 $$
where $\nu_j$ and $\eta_j$ are given by \eqref{e.nu} and
\eqref{e.eta}.

\cle\label{le:dist-int} If the first $(\sigma,\delta)$-hyperbolic
time map $h:M\to\ZZ^+$ belongs to $L^p(m)$ for some $p>4$, then
$\log\dist(x,\cs)$ is $\mu$-integrable. \fle

\dem We take any $\vare>0$ and use Proposition~\ref{p.dens} to
ensure the existence of two non-negative measures $\omega$ and
$\rho$ bounding $\eta_n$, where $\omega$ has  density bounded by
some constant and $\rho$ has total mass bounded by $\vare$. Recall
that $\rho$ was defined in \eqref{e.ro} by
 $$
 \rho=C_2\sum_{k=\ell}^\infty\sum_{j=1}^{k-1}
 f^j_*(m\mid H^*_k),
 $$
 where $\ell$ is some large integer.

Let us compute now  the weight  $\rho$ gives to some special
family of neighborhoods of $\cs$. For $i\ge 1$ let
$d_i=\sigma^{bi}$ where $0<\sigma<1$ comes from the definition of
$(\sigma,\delta)$-hyperbolic time. Define for $i\ge 1$
 $$
 B_i=\{x\in M\colon \dist(x,\cs)\le d_i\}.
 $$
 If $n$ is a
$(\sigma,\delta)$-hyperbolic time for $x\in M$, then $f^j(x)\in
M\setminus B_i$ for all $j\in\{n-i,\dots,n-1\}$. This implies that
\begin{align*}
\rho(B_i) & =  C_2 \sum_{k=\ell}^\infty\sum_{j=1}^{k-1} m(
H_k^*\cap f^{-j} (B_i ) )\\ & = C_2
\sum_{k=\ell}^\infty\sum_{j=1}^{k-i} m( H_k^*\cap f^{-j} (B_i ))
\\
& \le C_2 \sum_{k=\max\{\ell,i\}}^\infty\sum_{j=1}^{k-i} m(
H_k^*\cap f^{-j} (B_i) ) \\&\le C_2  \sum_{k=i}^\infty k \, m(
H^*_k),
\end{align*}
for all $i\ge1$. 
Now by Proposition~\ref{p.density} and Proposition~\ref{p.dens} we
know that
\[
\mu_n\le\frac1n\sum_{j=0}^{n-1}\nu_j+\omega+\rho\le\nu+\rho
\]
where $\nu$ is a measure with uniformly bounded density. Hence any
weak$^*$ accumulation point $\mu$ of the sequence $(\mu_n)_n$ is
bounded by $\nu+\rho$. Since we are assuming that $\cs$ is a
submanifold of $M$, then $\log\dist(x,\cs)$ is integrable with
respect to $\nu$ by Proposition~\ref{le:logdist-int}. On the other
hand,
 \begin{align*}
 \int_M -\log\dist_\delta(x,\cs)\,d\rho &\le
\sum_{i=1}^\infty -\rho(B_i)\log d_{i+1}  \le -b \log\sigma
\sum_{i=1}^\infty (i+1) \sum_{k=i}^\infty k \, m( H^*_k).
\end{align*}
We have $h\in L^p(m)$  by assumption, which is equivalent to
$\sum_{k\ge1} k^p m(H_k^*)<\infty.$ This implies that there is
some constant $C>0$ such that $m(H_k^*)\le C k^{-p}$ for all
$k\ge1$. Thus we have for $i\ge 2$
\[
\sum_{k=i}^\infty k \, m( H^*_k) \le \sum_{k=i}^\infty
\frac{C}{k^{p-1}}\le \int_{i-1}^\infty \frac{C}{x^{p-1}} \, dx =
\frac{C}{(p-2)(i-1)^{p-2}},
\]
and so
\[
\sum_{i=2}^\infty (i+1) \sum_{k=i}^\infty k \, m( H^*_k) \le
\frac{C}{p-2}\sum_{i=2}^\infty\frac{i+1}{(i-1)^{p-2}}.
\]
This last quantity is finite whenever $p>4$. Hence
$\log\dist(x,\cs)$ is integrable with respect to $\mu$ for all
$p>4$.
 \cqd

 As a consequence of the last results and the assumption that $f$
 behaves like a power of the distance near the exceptional set
 $\cs$, we obtain the result below.

\cco\label{co:dist-int} If the first $(\sigma,\delta)$-hyperbolic
time map $h$ belongs to $L^p(m)$ for
some $p>4$, 
then $ \log \| Df(x)^{-1} \|$ is $\mu$-integrable. \fco

\dem It is an easy consequence of condition (s$_1$) in the
definition of non-degenerate exceptional set that for some
$\zeta>\beta$ we have
\[
\big| \log \| Df(x)^{-1} \| \big| \le \zeta \big| \log\dist(x,\cs) \big|
\]
for all $x$ in a small open neighborhood $V$ of $\cs$. Hence
\[
\int_V \big| \log \| Df(x)^{-1} \| \big| \,d\mu \le \zeta \int_V
-\log\dist(x,\cs) \, d\mu < \infty,
\]
and since $\log \| Df(x)^{-1} \|$ is bounded on the compact set
$M\setminus V$, this function is necessarily integrable with
respect to $\mu$  on $M$. \cqd

Now we are ready to conclude the proof of Theorem~\ref{t.integrable-h-2} and Theorem~\ref{co:t-integrable-h}.
Assuming that the first $(\sigma,\delta)$-hyperbolic time map $h$ belongs to $L^p(m)$ for some $p>4$, it follows
from Lemma~\ref{le:dist-int} and Corollary~\ref{co:dist-int} that both:
\begin{enumerate}
    \item $\log\dist(x,\cs)$ is integrable with respect to $\mu$;
    \item $\log \| Df(x)^{-1} \|$ is integrable with respect to
    $\mu$.
\end{enumerate}
Observe that by definition of $(\sigma,\delta)$-hyperbolic time,
if $n$ is a $(\sigma,\delta)$-hyperbolic time for $x$ and if $k$
is a $(\sigma,\delta)$-hyperbolic time for $f^n(x)$, then $n+k$ is
a $(\sigma,\delta)$-hyperbolic time for $x$. Moreover, since $h$
is well defined Lebesgue almost everywhere and $f$ preserves sets
of Lebesgue measure zero, then Lebesgue almost all points must
have infinitely many hyperbolic times. Thus we have
\begin{equation}\label{e.nova}
\liminf_{n\to+\infty} \frac1n\sum_{j=0}^{n-1} \log \| Df(f^j(
x))^{-1} \| \le \log\sigma <0
\end{equation}
for Lebesgue almost every $x\in M$, and hence for $\mu$ almost
every $x\in M$.  The $\mu$-integrability of $\log \| Df(x)^{-1}
\|$ and Birkhoff's ergodic theorem then ensure that
\begin{equation}\label{eq:limit<0}
 \lim_{n\to+\infty}\frac1n\sum_{j=0}^{n-1} \log \| Df(f^j (x))^{-1}
\| \le \log\sigma<0
\end{equation}
for $\mu$ almost every $x\in M$.

\cre\label{re:corollary} If $\cs$ is equal to empty set, then $\log \| Df(x)^{-1} \|$ is immediately integrable
with respect to $\mu$ because it is a bounded function. Hence \eqref{eq:limit<0} is enough for obtaining
Theorem~\ref{co:t-integrable-h} by applying Theorem~\ref{t.abv2}. \fre

The strong integrability condition on $h$ ensures that
there exists an absolutely continuous invariant measure
$\mu$ and that $\log\dist(x,\cs)$ is $\mu$-integrable after
Lemma~\ref{co:dist-int}. Then we are in the setting of the
following result.

\cle\label{le:both-averages} If $\mu$ is an $f$-invariant probability measure and $\log\dist(x,\cs)$ is
$\mu$-integrable, then for every $0<\eta<1$ there is a set $R$ with $\mu(R)>1-\eta$ such that points in  $R$
have slow recurrence to $\cs$. \fle

\dem
We start by fixing a small $\eta>0$ and choosing $\alpha>0$ such
that
\[
\prod_{n\ge1} (1-e^{-\alpha n})\ge 1-\eta.
\]
The integrability of
$\log\dist(x,\cc)$ with respect to $\mu$
and the definition of the $\delta$-truncated distance
$\dist_{\delta}$ ensure that for each $k\in\NN$ we may find
$\delta_k>0$ for which
\[
\int_M -\log \dist_{\delta_k}(x,\cc) \, d\mu \le \frac1{k2^{k+1}}.
\]
We define for each $k\in\NN$
 $$
 \varphi_k(x)=\lim_{n\to+\infty}\frac1n\sum_{j=0}^{n-1}-\log
 \dist_{\delta_k}(f^j(x),\cc).
 $$
This $\varphi_k$ is well-defined $\mu$ almost everywhere in $M$ by
 Birkhoff's ergodic theorem. Moreover
 $$
 \int_M\varphi_k \,d\mu=\int_M-\log \dist_{\delta_k}(x,\cc)\,
 d\mu \le \frac1{k2^{k+1}}.
 $$
 Let
$$
E_k=\left\{ x\in M : \varphi_k(x)>\frac{1}{k}\right\}.
$$
Since $\varphi_k\ge 0$  we have
\[
\frac{\mu(E_k)}{k}\le \int_{E_k} \varphi_k \,d\mu \le
\int_{M}\varphi_k \,d\mu \le \frac1{k2^{k+1}},
\]
which implies that $\mu(E_k)\le 2^{-(k+1)}$. Hence we may find a big enough $k_1\in\NN$ such that
$\mu(M\setminus E_{k_1})\ge 1-e^{-\alpha}>0$. This is the first step in the following construction by induction
on $n$. Assuming that we have chosen $k_1<k_2<\dots<k_n$ satisfying
\[
\mu(M\setminus(E_{k_1}\cup\dots\cup E_{k_j})) \ge (1-e^{-\alpha j}) \mu(M\setminus( E_{k_1}\cup \dots\cup
E_{k_{j-1}}))>0
\]
for all $j=2,\dots,n$, then we may find a big enough
$k_{n+1}>k_n$ such that
\[
\mu(M\setminus(E_{k_1}\cup \dots\cup E_{k_n}\cup E_{k_{n+1}})) \ge (1-e^{-\alpha(n+1)}) \mu(M\setminus
(E_{k_1}\cup \dots\cup E_{k_n}))>0.
\]
Now, taking \( R=M\setminus\cup_{k\ge1} E_{k_n} \) we have \[ \mu(R) \ge \prod_{n\ge1} (1-e^{-\alpha n}) \ge
1-\eta.
\]
Let us now show that points in $R$ have slow approximation to $\cc$. Given $\vare>0$ we choose $n\in\NN$ for
which $\vare>1/k_n$. If $x\in R$, then in particular $x\notin E_{k_n}$, and this implies
 \[
 \lim_{n\to+\infty}\frac1n\sum_{j=0}^{n-1}-\log
 \dist_{\delta_{k_n}}(f^j(x),\cc)=
 \varphi_{k_n}(x)\le \frac1{k_n} < \vare  .
 \]
This concludes the proof of the result. \cqd

Since $\mu$ is absolutely continuous with respect to $m$, we deduce from \eqref{eq:limit<0} and
Lemma~\ref{le:both-averages} that there is a set with positive Lebesgue measure on which $f$ is non-uniformly
expanding  and whose points have slow recurrence to $\cs$.
 Actually these conditions must hold for
Lebesgue almost every $x\in M$. Indeed, let $H$ be the set of points for which both \eqref{liminf1} and
\eqref{e.faraway1} hold, and take $B=M\setminus H$. Observe that $B$ is invariant by $f$ and $h\in L^p(m\vert
B)$ with $p>4$. If $m(B)>0$, then  by the previous arguments we would prove the existence of some $A\subset B$
with $m(A)>0$ where $f$   is non-uniformly expanding  and points have slow recurrence to $\cs$. This would
naturally give a contradiction.

Thus we have proved that $f$ is non-uniformly expanding and points
have slow recurrence to $\cs$ Lebesgue almost everywhere. Applying
Theorem~\ref{t.abv2} we prove Theorem~\ref{t.integrable-h-2}.

\cre The hypothesis of $h$ belonging to $L^p(m)$ for some  $p>4$
can be replaced by $d\mu/dm\in L^q(m)$ for some $q>1$. In fact,
the  integrability of $h$ has only been used to prove that
 $\log\dist(x,\cc)$ is integrable with respect to $\mu$ --- which
 implies that
    $\log \| Df(x)^{-1} \|$ is also integrable with respect to
    $\mu$ by Remark~\ref{co:dist-int}.
As stated in
 Corollary~\ref{c.integra} this is a consequence of $d\mu/dm\in L^q(m)$
for some $q>1$.

\fre


\section{An example with non-integrable first hyperbolic time map}
\label{sec.example}

\begin{figure}[htbp]
  \centering
  \includegraphics[width=7.5cm, height=7cm]{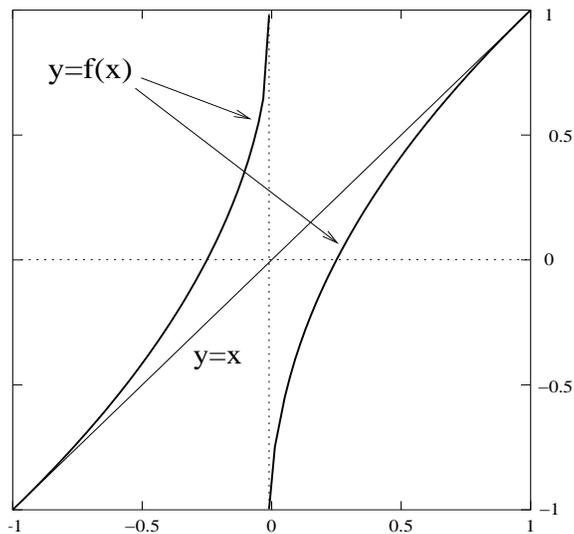}
  \caption{A map with non-integrable first hyperbolic
  time map.}
  \label{fig1}
\end{figure}

In this section we exhibit a map of the circle, differentiable
everywhere except at a single point, having a positive frequency
of hyperbolic times at Lebesgue almost every point, but whose
first hyperbolic time map is not integrable with respect to the
Lebesgue measure. This example is an adaptation of the
``intermittent'' Manneville map 
into a local homeomorphism of the circle; see~\cite{mp}. Consider
$I=[-1,1]$ and the map $\hat f:I\to I$ (see figure~\ref{fig1})
given by
\[
x\mapsto
  \begin{cases}
    2\sqrt{x}-1 & \text{if $x\ge0$}, \\
   1- 2\sqrt{|x|} & \text{otherwise}.
  \end{cases}
\]
This map induces a continuous local homeomorphism $f:S^1\to S^1$
through the identification $S^1 = I / \sim$, where $-1\sim 1$, not
differentiable at the point $0$.

\smallskip

\subsubsection*{Topological mixing} We will show that  given any open interval $J\subset S^1$ there is some
$N\in\NN$ such that $f^N(J)=S^1$. Note that $f$ has two inverse branches $ g_1:(-1,1)\to(0,1)$ and $
g_2:(-1,1)\to(-1,0)$, given by
$$
g_1(x)=\left(\frac{1+x}2\right)^2 \quad\mbox{and}\quad g_2(x)=-\left(\frac{1-x}2\right)^2.
$$Let $X=\{g_1^n(0),
g_2^n(0)\colon n\ge0\}$ be a set of points in the pre-orbit of $1\in S^1$ and let $\emptyset\neq J\subseteq S^1$
be an open interval. Observe that if $X\cap J\neq\emptyset$, then there is $n\ge1$ such that $1\in f^n(J)$, thus
the interval $f^n(J)$ would contain a neighborhood of $1$ in $S^1$ . We easily see that this implies
$f^{n+k}(J)=S^1$ for some $k\in\NN$. Hence to prove topological mixing for $f$ it is enough to show that given
any open interval $J\subset S^1$ there is $j\in\NN$ such that $f^j(J)\cap X\neq\emptyset$.


Let us take an interval $J\subset S^1$ such that $J\cap X=\emptyset$. Hence $0\not\in J$. Assume for
definiteness that $J\subset (-1,0)$. Thus there is $n\in\NN$ such that $f^n\vert_J$ is a diffeomorphism and
$f^n(J)\subset (0,1)$. It is clear that there is $\sigma>1$ independent of $n$ such that $m(f^{n+1}(J))\ge\sigma
m(J)$. Now let $J_1=f^{n+1}(J)\subseteq(0,1)$. If $J_1\cap X\neq\emptyset$, then we are done. Otherwise, by the
symmetry of $f$, we repeat the argument obtaining an iterate $J_2\subset(-1,0)$ of $J$ with
$m(J_2)\ge\sigma^{2}m(J)$. Since $(\sigma^{k})_{k\ge1}$ is unbounded, after a finite number of iterates the
image of $J$ will eventually hit $X$.

\smallskip

\subsubsection*{Invariance of Lebesgue measure}
Now we show that   Lebesgue measure is invariant by $f$. Note  that $f'(x)=|x|^{-1/2}$ for $x\in
S^1\setminus\{0\}$. It is straightforward to check that for all $x\in S^1\setminus\{0\}$
\[
\frac1{f'( g_1(x) )} + \frac1{ f'( g_2(x) )} = \sqrt{\left(\frac{1+x}2\right)^2} +
\sqrt{\left(\frac{1-x}2\right)^2} = 1,
\]
where \( g_1:(-1,1)\to(0,1)\) and \( g_2:(-1,1)\to(-1,0)\) are the inverse branches of $f$.
 Hence, the transfer operator
$$T_f: L^1(m)\to L^1(m)$$ given by
\[
T_f (\varphi)(x) = \sum_{f(z)=x} \frac{\varphi(z)}{|f'(z)|} =
\frac{\varphi(g_1(x))}{|f'(g_1(x))|} +
\frac{\varphi(g_2(x))}{|f'(g_2(x))|}
\]
fixes the constant densities. This means that Lebesgue measure is $f$-invariant.

\subsubsection*{Positive
  frequency of hyperbolic times}
Observe that $\log | (f'(x))^{-1} |=\log \sqrt{|x|}$ is  Lebesgue integrable on $S^1$ and that
\[
\int_{S^1} \log | (f'(x))^{-1} | \,dm = \int_{-1}^1 \frac12\log|x|
\,\left(\frac12 dx\right)= \frac12\int_0^1\log x\,dx=-\frac12
\]
(recall that we are taking the Lebesgue measure $m$ normalized on $S^1$). Thus, the invariance of Lebesgue
measure and Birkhoff's Ergodic Theorem ensure that
\[
G=\left\{x\in S^1 : \lim_{n\to+\infty}\frac1n\sum_{j=0}^{n-1}\log | (f'(f^j(x)))^{-1} | \le -\frac1{10} \right\}
\]
has positive Lebesgue measure. On the other hand, taking $\cs=\{0,\pm 1\}$ (the set of points where $f$ fails to
be a $C^2$ local diffeomorphism), we have
\begin{equation}
  \label{eq:slow-approx}
\int_{S^1} -\log\dist(x,\cs)\, dm(x) =2\int_0^{1/2} -\log x \,dx,
\end{equation}
and so $\log\dist(x,\cs)$ is  integrable with respect to $m$ on $S^1$.  Since $m$ is $f$-invariant,
Lemma~\ref{le:both-averages} ensures that there exists a set of points $R$ with slow recurrence to $\cs$ whose
complement has arbitrarily small Lebesgue measure. We obtain a  positively invariant subset $H=G\cap R$ with
$m(H)>0$ whose points satisfy conditions \eqref{liminf1} and \eqref{e.faraway1}. Then, by Lemma~\ref{l.disco} we
have that there is an interval $J\subset H$, up to a null Lebesgue measure subset. Due to the topological mixing
property and the regularity  of $f$ (preserves null Lebesgue measure sets) this implies that conditions
\eqref{liminf1} and \eqref{e.faraway1} hold Lebesgue almost everywhere on $S^1$. Hence, using
Theorem~\ref{t.abv2} one concludes that $f$ has positive frequency of hyperbolic times on Lebesgue almost every
point.


\subsubsection*{Non-integrability of the first hyperbolic time map}
We now show that given $0<\sigma<1$ and $\delta>0$  the  map
$h:S^1\to\ZZ_+$ assigning to each $x\in S^1$ the first
$(\sigma,\delta)$-hyperbolic time of $x$ cannot be Lebesgue
integrable in $S^1$. We observe that for $n$ being a
$(\sigma,\delta)$-hyperbolic time for $x$ it must satisfy
$$|f'(f^{n-1}(x))|\ge \sigma^{-1} > 1.$$ Hence  the first
$(\sigma,\delta)$-hyperbolic time for a given $x\in S^1\setminus\cs$ is at least the number of iterates  $x$
needs to hit a fixed neighborhood of $0$.

If we consider the inverse branch $g_1$ of $f$ and iterate a point $x_1\in(0,1)$ under $g_1$, we obtain a
sequence $(x_n)_{n\ge1}$ in $(0,1)$ satisfying
 \begin{equation}\label{e.rec}
 x_{n+1}=\frac{(1+x_n)^2}4,\quad n\ge1.
 \end{equation}
According to the observation above, we must have
\[
\int_{S^1} h\, dm \gtrsim \sum_{n\ge1} n (x_{n+1}-x_n).
\]
The non-integrability of $h$ is  then a consequence of the
following result.

\cle   ${\sum_{n\ge1} n (x_{n+1}-x_n)=+\infty}$.
  \fle
  \dem
We first prove (by induction) that
 \begin{equation}\label{e.xn}
0\le x_n\le 1-\frac1{2n}\quad\text{for every $n\ge 1$}.
 \end{equation}
 This obviously holds for $n=1$ since we have chosen $x_1\in (0,1/2)$. Assuming that
\eqref{e.xn} holds for $n\ge1$ we then have
\begin{align*}
0\le x_{n+1}&=\frac{(1+x_n)^2}4\\&\le \frac{(2-1/(2n))^2}4 \\&=
 1-\frac1{2n}+\frac{1}{16n^2}\\& = 1-\frac1{2n+2}
 \left(\frac{n+1}{n}-\frac{n+1}{8n^2}\right).
\end{align*}
 It is enough to observe that
 $$
 \frac{n+1}n-\frac{n+1}{8n^2}=\frac{8n^2+7n-1}{8n^2}>1,\quad\text{for all $n\ge1$.}
 $$
  Using the recurrence relation \eqref{e.rec}, a simple calculation now shows that
 $$x_{n+1}-x_n=\frac{(1-x_n)^2}4,$$
 which together with \eqref{e.xn} leads to
 $$
 x_{n+1}-x_n\ge\frac1{16n^2}.
 $$
 This is enough for concluding the proof of the result.
  \cqd




\begin{thebibliography}{ABV}

\bibitem{Al} J. F. Alves, {\em SRB measures for
non-hyperbolic systems with multidimensional expansion}, Ann.
Scient. Éc. Norm. Sup., $4^e$ série, {\bf 33} (2000), 1-32.

\bibitem{AA} J. F. Alves, V. Araújo, {\em Random perturbations
of non-uniformly expanding maps}, Astérisque 286 (2003), 25-62.


\bibitem{ABV} J. F. Alves, C. Bonatti, M. Viana, {\em
SRB measures for partially hyperbolic systems with mostly
expanding central direction}, Invent. Math. {\bf 140} (2000),
351-398.

\bibitem{ALP} J. F. Alves, S. Luzzatto, V. Pinheiro, {\em Markov
structures and decay of correlations for non-uniformly expanding dynamical systems}, preprint CMUP 2002.

\bibitem{AV} J. F. Alves, M. Viana, {\em Statistical
stability for robust classes of maps with nonuniform expansion},
Ergod. Th. \& Dynam. Sys. {\bf 22} (2002), 1-32.



\bibitem{BBM} V. Baladi, M.
Benedicks, V. Maume-Deschamps, {\em Almost sure rates of mixing
for i.i.d. unimodal maps}, Ann. Scient. Éc. Norm. Sup., $4^e$
série, {\bf 35}  (2002), 77-126.

\bibitem{BC1} M. Benedicks, L. Carleson, {\em On
iterations of $1-ax^2$ on $(-1,1)$}, Ann. Math. {\bf 122} (1985),
1-25.

\bibitem{BC2} M. Benedicks, L. Carleson, {\em The
dynamics of the H\'enon map}, Ann. Math. {\bf 133} (1991), 73-169.

\bibitem{BY1} M. Benedicks, L.-S. Young, {\em
Absolutely continuous invariant measures and random perturbations for certain one-dimensional maps}, Erg. Th. \&
Dyn. Sys. {\bf 12} (1992), 13-37.

\bibitem{BY2} M. Benedicks, L.-S. Young, {\em
SRB-measures for certain H\'enon maps}, Invent. Math. {\bf 112} (1993), 541-576.

\bibitem{DKK} D. Dolgopyat, V. Kaloshin, L. Koralov, {\em Hausdorff dimension in stochastic dispersion}, J.
Statyst. Phys. {\bf 108}, No. 5-6 (2002) 943-971.

\bibitem{J} M. Jakobson, {\em  Absolutely continuous
invariant measures for one-parameter families of one-dimensional
maps}, Comm. Math. Phys. {\bf 81}
 (1981), 39-88.

\bibitem{KH}
A. Katok, B. Hasselblatt,
\newblock {\em Introduction to the modern theory of dynamical
systems},
\newblock Cambridge Univ. Press, New York, 1995.


\bibitem{Man87}
R.~Ma{\~n}{\'e}.
\newblock {\em Ergodic theory and differentiable dynamics}.
\newblock Springer Verlag, 1987.

\bibitem{mp} P. Manneville,
\emph{Intermittency, self-similarity and $1/f$ spectrum in
dissipative  dynamical systems}, J. Physique \textbf{41} (1980), 1235-1243.

\bibitem{MV} L. Mora, M. Viana, {\em Abundance of
strange attractors}, Acta Math. {\bf 171} (1993), 1-71.


\bibitem{Pli} V. Pliss, {\em On a conjecture due to
Smale}, Diff. Uravnenija, {\bf 8} (1972), 262-268.

\bibitem{V1} M. Viana, {\em Strange attractors in
higher dimensions}, Bull. Braz. Math. Soc. {\bf 24} (1993), 13-62.

\bibitem{V} M. Viana, {\em Multidimensional
non-hyperbolic attractors}, Publ. Math. IHES {\bf 85} (1997),
63-96.

\bibitem{WY} Q. Wang, L.-S. Young, {\em Strange attractors with one direction of instability}, Commun. Math.
Phys. {\bf 218}, No.1 (2001) 1-97.


\end{thebibliography}
\end{document}